\documentclass[11pt]{eptcs}

\usepackage{xcolor}
\definecolor{bluegreen}{rgb}{0.0, 0.62, 0.24}

\usepackage{environ}

\makeatletter
\newsavebox{\measure@tikzpicture}
\NewEnviron{scaletikzpicturetowidth}[1]{%
	\def\tikz@width{#1}%
	\def\tikzscale{1}\begin{lrbox}{\measure@tikzpicture}%
		\BODY
	\end{lrbox}%
	\pgfmathparse{#1/\wd\measure@tikzpicture}%
	\edef\tikzscale{\pgfmathresult}%
	\BODY 
}
\makeatother

\usepackage[centertags,sumlimits,intlimits,namelimits,reqno]{amsmath}
\usepackage{amsfonts,amsthm,amssymb}
\usepackage{mathtools}
\usepackage{enumitem}
\usepackage{dsfont}

\usepackage{tikz}
\usepackage{pgfplots}
\pgfplotsset{compat=1.15}
\usepgfplotslibrary{fillbetween}
\usetikzlibrary{patterns}
\usetikzlibrary{backgrounds}
\usetikzlibrary{arrows.meta,decorations.markings}

\usepackage{float}
\usepackage{wrapfig}

\usepackage[capitalize]{cleveref}
\usepackage[backend=bibtex, bibencoding=utf8, maxbibnames = 10, style = alphabetic]{biblatex}

\usepackage{subfig}

\hypersetup{colorlinks, linkcolor=black!50!blue, citecolor=black!50!blue, urlcolor=black!50!blue, breaklinks=true}

 \usetikzlibrary{ 
  	cd,
  	math,
  	decorations.markings,
	decorations.pathreplacing,
	decorations.pathmorphing,
  	positioning,
  	arrows.meta,
  	shapes,
	shadows,
	shadings,
	patterns,
  	calc,
  	fit,
  	quotes,
  	intersections,
  }
\tikzset{dotsF/.style={pattern=north east lines,pattern color=black!60!white}}

\tikzset{
	oriented WD/.style={
		every to/.style={out=0,in=180,draw},
    label/.style={
    	font=\everymath\expandafter{\the\everymath\scriptstyle},
      inner sep=0pt,
      node distance=2pt and -2pt},
    semithick,
    node distance=1 and 1,
    decoration={markings, mark=at position \stringdecpos with \stringdec},
    ar/.style={postaction={decorate}},
    execute at begin picture={\tikzset{
    	x=\bbx, y=\bby,
      every fit/.style={inner xsep=\bbx, inner ysep=\bby}}}
    },
    string decoration/.store in=\stringdec,
    string decoration={\arrow{stealth};},
    string decoration pos/.store in=\stringdecpos,
    string decoration pos=.7,
    bbx/.store in=\bbx,
    bbx = 1.5cm,
    bby/.store in=\bby,
    bby = 1.5ex,
    bb port sep/.store in=\bbportsep,
    bb port sep=1.5,
    bb port length/.store in=\bbportlen,
    bb port length=4pt,
    bb penetrate/.store in=\bbpenetrate,
    bb penetrate=0,
    bb min width/.store in=\bbminwidth,
    bb min width=1cm,
    bb rounded corners/.store in=\bbcorners,
    bb rounded corners=2pt,
    bb small/.style={
    	bb port sep=1, bb port length=2.5pt, bbx=.4cm, bb min width=.4cm, bby=.7ex},
		bb medium/.style={
			bb port sep=1, bb port length=2.5pt, bbx=.4cm, bb min width=.4cm, bby=.9ex},
    bb/.code 2 args={
    	\pgfmathsetlengthmacro{\bbheight}{\bbportsep * (max(#1,#2)+1) * \bby}
      \pgfkeysalso{draw,minimum height=\bbheight,minimum
       width=\bbminwidth,outer sep=0pt,
         rounded corners=\bbcorners,thick,
         prefix after command={\pgfextra{\let\fixname\tikzlastnode}},
         append after command={\pgfextra{\draw
            \ifnum #1=0{} \else foreach \i in {1,...,#1} {
            	($(\fixname.north west)!{\i/(#1+1)}!(\fixname.south west)$) +(-\bbportlen,0) coordinate (\fixname_in\i) -- +(\bbpenetrate,0) coordinate (\fixname_in\i')}\fi 
            \ifnum #2=0{} \else foreach \i in {1,...,#2} {
            	($(\fixname.north east)!{\i/(#2+1)}!(\fixname.south east)$) +(-
\bbpenetrate,0) coordinate (\fixname_out\i') -- +(\bbportlen,0) coordinate (\fixname_out\i)}\fi;
           }}}
		},
			bb name/.style={
     	append after command={
				\pgfextra{\node[anchor=north] at (\fixname.north) {#1};}
			}
		}
  }

\tikzset{
	unoriented WD/.style={
  	every to/.style={draw},
  	shorten <=-\penetration, shorten >=-\penetration,
  	label distance=-2pt,
  	thick,
  	node distance=\spacing,
  	execute at begin picture={\tikzset{
  		x=\spacing, y=\spacing}
		}
  },
  pack size/.store in=\psize,
  pack size = 8pt,
	penetration/.store in=\penetration,
	penetration = 2pt,
  spacing/.store in=\spacing,
  spacing = 8pt,
  link size/.store in=\lsize,
  link size = 2pt,
  pack color/.store in=\pcolor,
  pack color = blue,
 	pack inside color/.store in=\picolor,
  pack inside color=blue!20,
 	pack outside color/.store in=\pocolor,
  pack outside color=blue!50!black,
 	surround sep/.store in=\ssep,
  surround sep=8pt,
 	link/.style={
  	circle, 
  	draw=black, 
  	fill=black,
  	inner sep=0pt, 
 		minimum size=\lsize
 	},
  pack/.style={
 		circle, 
 		draw = \pocolor, 
  	fill = \picolor,
  	inner sep = .25*\psize,
 		minimum size = \psize
  },
  outer pack/.style={
 		ellipse, 
 		draw,
  	inner sep=\ssep,
  	color=\pocolor,
 	},
  intermediate pack/.style={
 		ellipse,
 		dashed, 
  	draw,
  	inner sep=\ssep,
 		color=\pocolor,
 	},
}
  
\tikzset{
	spider diagram/.style={
		every to/.style={out=0, in=180, draw, thick},
		thick
	},
	dot size/.store in=\dotsize,
	dot size = 5pt,
	dot fill/.store in=\dotfill,
	dot fill = black,
	leg length/.store in=\leglen,
	leg length = 15pt,
	baby/.style={dot size = 2pt, leg length = 6pt},
	young/.style={dot size = 3pt, leg length = 10pt},
	special spider/.code n args={4}{
		\pgfkeysalso{circle, draw, thick, inner sep=0, fill=\dotfill, minimum width=\dotsize,
  		prefix after command={\pgfextra{\let\fixname\tikzlastnode}},
  		append after command={\pgfextra{
  			\ifnum #1=0{} \else {\foreach \i in {1,...,#1} {
					\tikzmath{\anglei={-90*(#1+1-2*\i)/#1};}
  				\draw [thick]
						(\fixname) .. controls 
						($(\fixname.center)-(\anglei:#3/3)$) and ($(\fixname.center)-(\anglei:#3*2/3)$) .. 
						({$(\fixname.center)-(\anglei:#3*2/3)$}-|{$(\fixname.center)-(#3,0)$}) coordinate (\fixname_in\i);
  			}}\fi
  			\ifnum #2=0{} \else {\foreach \i in {1,...,#2} {
					\tikzmath{\anglei={90*(#2+1-2*\i)/#2};}
  				\draw [thick]
						(\fixname.center) .. controls 
						($(\fixname.center)+(\anglei:#4/3)$) and ($(\fixname.center)+(\anglei:#4*2/3)$) .. 
						({$(\fixname.center)+(\anglei:#4*2/3)$}-|{$(\fixname.center)+(#4,0)$}) coordinate (\fixname_out\i);
  			}}\fi
  		}}
		}
	},
	spider/.code 2 args={
		\pgfkeysalso{special spider={#1}{#2}{\leglen}{\leglen}}
	}
}

\tikzset{
  function/.style={->, thin, shorten <=4pt, shorten >=4pt}
}

\tikzset{
  tick/.style={
  	postaction={
    	decorate,
      decoration={
      	markings, mark=at position 0.5 with {
					\draw[-] (0,.4ex) -- (0,-.4ex);
				}
			}
		}
	}
}

\setlist{itemsep=-1pt}
\setlist[description]{leftmargin=!, itemindent=0em, parsep=0em, topsep=0em, listparindent=2em, labelwidth=2em}

\newtheorem*{theorem*}{Theorem}
\newtheorem{definition}{Definition}[section]

\theoremstyle{remark}

\newtheorem*{example*}{Example}
\newtheorem{remark}[definition]{Remark}

\crefalias{chapter}{section}

\title{Temporal Landscapes: A Graphical Logic of Behavior\thanks{B. Fong and D.I. Spivak were supported by AFOSR grant FA9550--17--1--0058.}}

\author{Brendan Fong\thanks{Email: brendan@topos.institute}\institute{Topos Institute} \and Alberto Speranzon\thanks{Email: alberto.speranzon@gmail.com} \institute{Honeywell Aerospace} \and David I.\ Spivak\thanks{Email:david@topos.institute}\institute{Topos Institute \& MIT, LIDS}}

\def\titlerunning{Temporal Landscapes: A Graphical Logic of Behavior}
\def\authorrunning{B. Fong, A. Speranzon \& D.I. Spivak}

\date{\vspace{-.3in}}

\usepackage{xparse}

\NewDocumentCommand{\True}{O{gray!40} O{1} m m}{
	\fill[draw=black, fill=#1, opacity=#2] (#3,#3) -- (#3,#4) -- (#4,#4);
}

\NewDocumentCommand{\AlwaysTrue}{O{gray} O{45} m m}{
	\fill[draw=white, bottom color = #1!40, top color = #1!0, middle color=#1!20, shading angle=#2] (#3,#3) -- (#3,#4) -- (#4,#4);
}

\newcommand{\LandscapeTimeAxis}[2]{
	\coordinate (ll) at (#1-0.5,#1-0.5);
	\coordinate (ur) at (#2+0.5, #2+0.5);	
	\draw[<->] (ll) -- (ur);
	\foreach \i in {0,..., #2}
		\draw (\i+.1,\i-.1) to[pos=-1.5, font=\tiny] node {\i} (\i-.1, \i+.1);
}

\NewDocumentCommand{\AlwaysTrueFrom}{O{gray} O{45} m m}{
	\fill[draw=white, bottom color = #1!40, top color = #1!0, middle color=#1!20, shading angle=#2] (#3,#3) -- (#3,#4) -- (#4,#4);
	\draw[draw=#1] (#3,#3) -- (#3,#4-1.5) edge[dashed] (#3,#4);
}

\newcommand{\LandscapeTimeAxisNoTicks}[2]{
	\coordinate (ll) at (#1-0.5,#1-0.5);
	\coordinate (ur) at (#2+0.5, #2+0.5);	
	\draw[<->] (ll) -- (ur);
}

\NewDocumentCommand{\False}{O{gray!40} O{1} m m}{
	\draw[ultra thick, double distance=0.03pt, draw=#1, opacity=#2] (#3,#3) -- (#4,#4);
}

\NewDocumentCommand{\AlwaysFalse}{O{gray!40} m m m}{
	\draw[ultra thick, dashed, double distance=0.03pt, draw=#1] (#2-#4,#2-#4) -- (#2,#2);
	\draw[ultra thick, double distance=0.03pt, draw=#1] (#2,#2) -- (#3,#3);
	\draw[ultra thick, dashed, double distance=0.03pt, draw=#1] (#3,#3) -- (#3+#4,#3+#4);
}

\NewDocumentCommand{\UnionAlwaysTrue}{ O{gray!40} O{1} m m m}{
	\fill[draw=black, fill=#1, opacity=#2] (#3,#3) -- (#3,#3+#5) -- (#4-#5,#4) -- (#4,#4);
}

\NewDocumentCommand{\ArbitraryLandscape}{O{gray!40} O{1} m}{
	\fill[draw=black, fill=#1, opacity=#2] #3;
}

\newcommand{\const}[1]{\mathtt{#1}}
\newcommand{\Cat}[1]{{\mathsf{#1}}}

\newcommand{\interval}[1]{{\downarrow\kern-0.4em#1}}

\newcommand{\wt}[1]{\widetilde{#1}}
\renewcommand{\ss}{\subseteq}

\newcommand{\pred}[1]{\mathsf{#1}}
\newcommand{\type}[1]{\mathsf{#1}}

\newcommand{\until}{\mathbin{\mathsf{U}}}

\newcommand{\true}{\const{true}}
\newcommand{\false}{\const{false}}
\newcommand{\Time}{\const{Time}}

\newcommand{\prop}{\const{Prop}}

\newcommand{\occ}{\pred{Occ}}
\newcommand{\free}{\pred{Free}}

\newcommand{\pos}{\pred{Pos}}
\newcommand{\velbnd}{\pred{v_{max}}}

\newcommand{\imp}{\Rightarrow}
\renewcommand{\iff}{\Leftrightarrow}

\newcommand{\smset}{\Cat{Set}}

\newcommand{\timebetw}[2]{\pred{TimeBetw}(#1,#2)}

\newcommand{\rr}{{\mathds{R}}}
\newcommand{\qq}{{\mathds{Q}}}

\newcommand{\nn}{\mathds{N}}
\newcommand{\zz}{\mathds{Z}}
\newcommand{\ir}{\mathds{I\hspace{.5pt}R}}

\newcommand{\commentout}[1]{}

\begin{document}   

\maketitle

\begin{abstract}
We present an elementary introduction to a new logic for reasoning about behaviors that occur over time. This logic is based on \emph{temporal type theory}. The syntax of the logic is similar to the usual first-order logic; what differs is the notion of \emph{truth value}. Instead of reasoning about whether formulas are true or false, our logic reasons about \emph{temporal landscapes}. A temporal landscape may be thought of as representing the set of durations over which a statement is true. To help understand the practical implications of this approach, we give a wide variety of examples where this logic is used to reason about autonomous agents.
\end{abstract}


\section{Introduction}

Logical formalization of temporal considerations has a long and rich history, though the first modern treatment is probably the tense logic of Prior \cite{prior1967past}, which has yielded what is now known simply as temporal logic. In temporal logic, as in any logical system, one has a syntactic way of building new formulas from simpler ones, say using conjunction and negation ($\varphi\wedge\neg\psi$), as well as a notion of model whose purpose is to specify the truth value of each such formula. 


Temporal logic gets its expressive power from various operators, which collect information about other times into the current time. For example, in linear temporal logic (LTL), one considers the binary \emph{until} operator $\until$. At time $t$, one may ask whether $\varphi$ will hold until $\psi$ holds, denoted $t\models\varphi\until \psi$, which more precisely means that there exists $t'>t$ such that $t'\models \psi$ and $t''\models\varphi$ for all $t<t''<t'$. One can understand all such operators as given by some sort of quantification over $t:T$, replacing each formula, e.g.\ $\varphi$, by a predicate $\varphi(t)$ in one variable. Restricting first-order logic by requiring that all atomic predicate symbols take only one variable, one obtains what is known as \emph{first-order monadic logic}, and adding the 2-ary predicate $t<t'$, one obtains what is known as the first order monadic logic \emph{of order}, $FO(<)$. It was shown by Kamp \cite{Kamp:1968a} that temporal logic with the until operator, together with its past-tense cousin \emph{since}, is precisely as expressive as $FO(<)$. Various additions and restrictions have been proposed over the years, in attempts to co-optimize between expressivity and computability. 

A completely different approach to temporal reasoning, known as \emph{temporal type theory} (TTT), was given in \cite{schultz2019temporal}. Instead of defining new logical operators that collate past and future times into the present, temporal type theory alters the very notion of truth itself, to make truth inherently depend on time. The goal of this article is to describe, in elementary terms, how TTT makes this idea precise, as well as how it can be used in practice. 

Temporal type theory begins by defining a topological space called the \emph{interval domain} $\ir$, whose points are the closed intervals $[t_1,t_2]\ss\rr$, which we call \emph{time-intervals}, and whose open sets are generated by open intervals $(a,b)$, each of which consists of all points $[t_1,t_2]$ with $a<t_2\leq t_2<b$. A \emph{sheaf} $B$ on $\ir$ is a type of behavior: it assigns to each basic open $(a,b)$ a set $B(a,b)\in\smset$, and for every $a\leq a'\leq b'\leq b$, it assigns a restriction function $\rho\colon B(a,b)\to B(a',b')$ that clips a longer-lasting behavior $x\in B(a,b)$ to a shorter-lasting behavior $\rho(x)\in B(a',b')$. Behavior types include:
\begin{enumerate}[topsep=0pt,itemsep=-1ex,partopsep=1ex,parsep=1ex]
	\item $\nn$, $\zz$, $\qq$, and $\rr$, the behaviors of natural numbers, integers, rationals, and reals (unchanging over any interval $(a,b)$);
	\item $\wt{\rr}$, the behavior of ``varying'' real numbers (changing continuously over any $(a,b)$);%
	\footnote{Note that the constant reals can be considered as a subtype $\rr\ss\wt{\rr}$ of the varying real numbers.}
	\item for any vector field $V$ on a topological space, the behavior of integral curves through $V$ (of duration $b-a$);
	\item for any graph $G$, the behavior of stochastically-timed walks through $G$;
	\item more generally, for any hybrid system \cite{henzinger2000theory}, the behavior of all legal trajectories;
	\item $\prop$, the behavior of truth values, also known as propositions, which one can think of as audits or monitors of behavior. $\prop$ will be the main character in this paper;
	\item the empty behavior and the singleton behavior ($\Time$ itself), as well as products, unions, subobjects, quotients, and exponentials of all the above.
\end{enumerate}
Discussing behavior types in detail is out of scope for this paper, as it includes definitions of sheaves and toposes; the interested reader is referred to \cite{schultz2019temporal} for a technical discussion, or to \cite[Chapter 7]{fong2019seven} for a gentle introduction. 
The goal in this paper is to give the reader a relatively self-contained understanding of $\prop$---the behavior type of truth values---in terms of \emph{temporal landscapes}.

We will also not give a detailed comparison between the expressive power of various temporal logics, such as Metric Interval Temporal Logic (MITL)~\cite{Alur.Feder.Henzinger:1996a} or Signal Temporal Logic (STL)~\cite{Maler.Nickovic:2004}, with that of temporal type theory. The main difference is simply that temporal type theory is a type theory, meaning that it can combine and reason about various types of behaviors, as exemplified above. Another is that temporal logic assumes a kind of omniscience about the future: the truth value of a proposition at time $t_1$ can contain information about what occurs over a whole interval $[t_1,t_2]$. TTT does not have this omniscience: a proposition whose truth value depends on more than one moment---such as ``whenever $A$ occurs at time $t_1$, $B$ must occur before time $t_2$''---is only falsifiable on long-enough intervals, e.g.\ those containing $[t_1,t_2]$. The aggregated truth value of a proposition over all intervals is its \emph{temporal landscape}; information about what occurs over an interval is ``stored'' over the interval itself, not at its left endpoint. However, LTL and MITL do embed as a fragment of TTT \cite[Chapter 8.6]{schultz2019temporal}, so proofs from these logics are valid in TTT.

Our focus will be on the \emph{descriptive power} of temporal type theory: through increasingly complex examples we shall demonstrate how TTT---in particular temporal landscapes---can be used to accurately model the relevant time-varying phenomena. We hope that this will be enough to give the reader a basic understanding of these ideas, even if translating them into reasoning power holds few subtleties.

Luckily, the base language of TTT is standard higher-order logic, which is quite similar to first-order logic. Not only should reasoning in TTT thus be more familiar to users with a grounding in predicate logic, a wide variety of proof assistants, including HOL, Lean, and Coq, can hence be easily adapted to provide formal verification of reasoning in TTT \cite{Coquand.Huet:1988a,deMoura.Soonho.Avigad.VanDoorn.vonRaumer:2015a,nipkow2002isabelle}.

We will begin in \cref{chap.temporal_landscapes} by defining temporal landscapes. In \cref{sec.logical_manipulations} we discuss logical operations on temporal landscape and in \cref{chap.examples} we give several increasingly expressive examples of temporal landscapes in the context of autonomous agents. Conclusions are provided in \cref{chap.conclusion}.

\section{Definition of temporal landscape}\label{chap.temporal_landscapes}

Temporal landscapes provide the \emph{truth values} of a logical system, which we call \emph{temporal landscape logic}. Truth values may be thought of as acceptable answers to ``yes/no''-style questions. For example, in standard propositional logic, the truth values are simply $\true$ and $\false$. In propositional logic then, the question ``Is it raining?'' may be answered with ``yes'' ($\true$) or ``no'' ($\false$). In temporal landscape logic, the answer to this question is a temporal landscape indicating precisely those time intervals during which it is raining. Let us be a bit more precise.

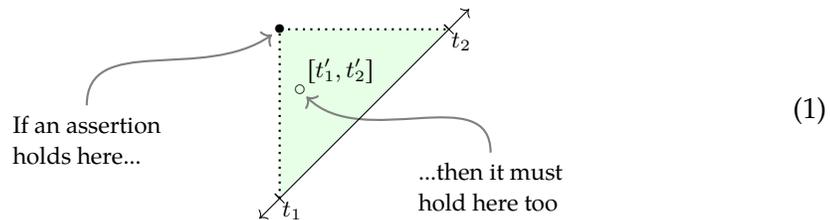
\begin{figure}[b!]
\centering
\vspace*{-12pt}
\subfloat[][]{\label{eqn.downclosure}
\begin{scaletikzpicturetowidth}{0.3\textwidth}
		\begin{tikzpicture}[scale=\tikzscale, baseline=(current  bounding  box.center)]
	\coordinate (ll) at (.7, .7);
	\coordinate (ur) at (3.8, 3.8);
	\coordinate (x) at (1,1);
	\coordinate (y) at (3.5,3.5);
	\coordinate (x') at (1.3,1.3);
	\coordinate (y') at (2.6,2.6);
	\draw[fill=green!10, dotted, thick] (x) -- (x|-y) -- (y);
	\draw (x) -- +(-2pt,2pt) to node[pos=2.8] {$t_1$} +(2pt,-2pt);
	\draw (y) -- +(-2pt,2pt) to node[pos=2.8] {$t_2$} +(2pt,-2pt);
	\node[pin={[pin edge={thick, <-,shorten <=-4pt, out=215, in=90}, align=left]215:{If an assertion\\ holds here...}}] at (x|-y) {$\bullet$};
	\node[label={[xshift=0.5cm, yshift=-0.2cm]${[t_1',t_2']}$}] (x'y') at (x'|-y') {$\circ$};
	\node[pin={[pin edge={thick, <-,shorten <=-2pt, out=315, in=90},align=left]330:{...then it must\\ hold here too}}] at (x'y') {};
	\draw[<->] (ll) -- (ur);
\end{tikzpicture}
\end{scaletikzpicturetowidth}}
\hspace*{1cm}
\subfloat[][]{\label{eqn.openness}
	\begin{scaletikzpicturetowidth}{0.3\textwidth}
		\begin{tikzpicture}[scale=\tikzscale, baseline=(current  bounding  box.center)]
		\coordinate (ll) at (.7, .7);
		\coordinate (ur) at (3.8, 3.8);
		\coordinate (x) at (1, 1);
		\coordinate (y) at (3.5, 3.5);
		\coordinate (x') at (1.6, 1.6);
		\coordinate (y') at (2.6, 2.6);
		\coordinate (if) at (1.4, 2.9);
		\coordinate (then) at ($(if)+(-0.25, 0.15)$);
		\draw[fill=green!10, dotted, thick] (x) -- (x|-y) -- (y);
		\node[pin={[pin edge={thick, <-,shorten <=-2pt, out=70, in=0},align=left, thick]175:{If an assertion\\ holds here...}}] at (if) {$\bullet$};
		\node[pin={[pin edge={thick, <-,shorten <=-4pt, out=215, in=0},align=left]200:{...then there exists\\ a larger interval\\ where it holds too}}] at (then) {$\circ$};
		\draw (x) -- +(-2pt,2pt) to node[pos=2.8] {$t_1$} +(2pt,-2pt);
		\draw (y) -- +(-2pt,2pt) to node[pos=2.8] {$t_2$} +(2pt,-2pt);
		\draw[<->] (ll) -- (ur);
	\end{tikzpicture}
	\end{scaletikzpicturetowidth}}
	\caption{Downclosure and openness.}
\end{figure}

Let $\rr$ be the real numbers, thought of as representing points in time. Given two times $t_1,t_2 \in \rr$ with $t_1\leq t_2$, we write $[t_1,t_2]$ for the set of all times between $t_1$ and $t_2$; we call this a \emph{time interval}. Graphically, we may represent a time interval by a point above the diagonal in the plane $\rr^2$, see~\cref{eqn.downclosure}.

A \emph{temporal landscape} is a set of time intervals with two special properties. The first is known as \emph{down-closure}: if a time interval $[t_1,t_2]$ is in the temporal landscape, and $[t_1',t_2']$ is contained in $[t_1,t_2]$, then $[t_1',t_2']$ is in the temporal landscape too. This property makes the assumption that if an assertion holds throughout a time interval, then it holds on all subintervals. For example, if it is raining throughout the time interval from 9:00 to 13:00, then it is also raining throughout the time interval from 10:00 to 10:45. In pictures, this means that a temporal landscape must be closed under both moving right and moving downward, as shown in \cref{eqn.downclosure}.

The second property of temporal landscapes is an \emph{openness} (sometimes called a \emph{roundedness}) property: if an assertion holds on some $[t_1,t_2]$, then there exists some larger interval, $[t'_1, t'_2]$ with \emph{both} $t_1'<t_1$ and $t_2<t_2'$. This larger interval may only be infinitesimally larger, but it must be strictly larger on both sides. In pictures, we illustrate this as in \cref{eqn.openness}.The above is summarized in~\cref{def.temp_landscape}.
\begin{definition}\label{def.temp_landscape}
	A \emph{temporal landscape on $\rr$} is a set $L$ of time intervals $[t_1,t_2] \subseteq \rr$, $t_1\leq t_2$, such that
	\begin{enumerate}[label=(\alph*)]
		\item if $[t_1,t_2] \in L$, and $t_1 \le t_1' \le t_2' \le t_2$, then $[t_1',t_2'] \in L$.
		\item if $[t_1, t_2]\in L$ then there exists $t_1' <t_1 \le t_2<t_2'$ such that $[t_1',t_2'] \in L$.
	\end{enumerate}
	We write $\prop$ for the set of temporal landscapes.
\end{definition}

Together, requirements (a) and (b) state that temporal landscapes form the open sets of the Scott topology on the \emph{interval domain} $\ir$, a well-studied topological space in domain theory \cite{Gierz.Keimel.Lawson.Mislove.Scott:2003a}. While we will not need any topos or sheaf theory here, we remark that sheaves on $\ir$ form a topos, whose subobject classifier consists precisely of temporal landscapes; this topos is the subject of \cite{schultz2019temporal}. 

\begin{remark}
	By definition---\cref{def.temp_landscape}(b)---a temporal landscape $L$ does not include its boundary: it is an open set in $\ir$. Hence in our examples so far, \cref{eqn.downclosure,eqn.openness}, we have drawn them using a dotted line. From now on we use the visually simpler convention of drawing them with a solid line.
\end{remark}

The simplest temporal landscape is that of a \emph{roof}; this form the basis for the topology on $\ir$.

\begin{definition}
Given a pair $a<b$ in $\rr$, the \emph{roof over $a,b$} is the temporal landscape 
\[
	\timebetw{a}{b} \coloneqq \{ [t_1,t_2] \mid a< t_1 \le t_2< b\}
\]
Given any pair of real numbers $a<b$, a \emph{temporal landscape on $(a,b)$} is a temporal landscape that is a subset of $\timebetw{a}{b}$.
\end{definition}

General temporal landscapes are curves that remain above the diagonal line and whose slope is piecewise continuous and remains in the interval $[0,\infty]$. Note that rotating by $45^\circ$, the slope condition becomes precisely the statement that the curve must be what is known as a 1-Lipschitz function.

\section{Temporal landscape logic}\label{sec.logical_manipulations}

Temporal landscapes form the truth values of a logical system. More precisely, temporal landscapes form the elements of what is known as a Heyting algebra. This means that standard logical constants and operations, such as $\true$, $\false$, AND ($\wedge$), OR ($\vee$) and implication ($\imp$), have interpretations as temporal landscapes and operations on them.

To begin, we introduce the temporal landscapes $\true$ and $\false$. The temporal landscape $\true$ contains all time intervals $	\true \coloneqq \{[t_1,t_2]\mid t_1 < t_2 \in \rr\}$. This landscape is the maximal one that visually can be depicted as a infinite triangle over the time line. 
On the other hand, the temporal landscape $\false$ contains no time intervals at all: $\false \coloneqq \varnothing$ and it is the minimal landscape.

Given temporal landscapes $\varphi$ and $\psi$, their conjunction $\varphi \wedge \psi$ is given by their intersection, and their disjunction $\varphi \vee \psi$ is given by their union. For an explicit example of conjunction,  see \cref{fig.grid_world_dynamic} on page \pageref{fig.grid_world_dynamic}. In that example, the temporal landscape $\free(\pred{Nbr}(v))$ on the right is the conjunction---that is, the intersection---of the temporal landscapes $\free(w_r)$ and $\free(w_u)$ on the left and center.

It is straightforward to check that the results of these operations are again temporal landscapes, and that $\wedge,\vee$ obey the usual properties of (constructive) first-order logic; for example, for any temporal landscape $\varphi$, we have $\varphi \wedge \true = \varphi$.

Defining an implication that obeys the usual properties is a bit more subtle. Given temporal landscapes $\varphi$ and $\psi$, we define the temporal landscape
\[
	(\varphi \imp \psi) \coloneqq \{[a,b] \mid \timebetw{a}{b} \cap \varphi \subseteq \psi\}.
\]
To become acquainted with implication in general, we start with a special case, namely that of negation.

\begin{figure}[t!]
	\centering
	\[
	\left(
	\begin{aligned}
		\begin{scaletikzpicturetowidth}{.21\textwidth}
			\begin{tikzpicture}[scale=\tikzscale]
				\True[blue!40][0.5]{3}{6}
				\True[blue!40][0.5]{7}{15}
				\LandscapeTimeAxis{0}{15}
			\end{tikzpicture}
		\end{scaletikzpicturetowidth}
	\end{aligned}
	\hspace{-.54in}\imp\hspace{-.2in}
	\begin{aligned}
		\begin{scaletikzpicturetowidth}{.21\textwidth}
			\begin{tikzpicture}[scale=\tikzscale]
				\UnionAlwaysTrue[red!40][0.5]{0}{15}{3}
				\LandscapeTimeAxis{0}{15}
			\end{tikzpicture}
		\end{scaletikzpicturetowidth}
	\end{aligned}
	\right)
	\quad = \quad
	\begin{aligned}
		\begin{scaletikzpicturetowidth}{.21\textwidth}
			\begin{tikzpicture}[scale=\tikzscale]
				\ArbitraryLandscape{(0,0) -- (0,10) -- (7,10) -- (12,15) -- (15,15)}
				\draw[dashed] (7,7) --(7,10);
				\draw[dashed] (7,10) -- (10,10);
				\LandscapeTimeAxis{0}{15}
			\end{tikzpicture}
		\end{scaletikzpicturetowidth}
	\end{aligned}
	\]
	\caption{Temporal landscape for the implication $\varphi \imp \psi$.}\label{fig.temp_land_implication}
\end{figure}

The negation operator $\neg \varphi$ is given by $\neg \varphi \coloneqq (\varphi \imp \false)$. Equivalently, since $\false$ is the empty set, we may write
\[
	\neg \varphi = \{[a,b] \mid \timebetw{a}{b} \cap \varphi = \varnothing\}.
\] 


The visual intuition of the implication $\varphi \imp \psi$ generalizes that of negation, replacing the empty landscape $\false$ with $\psi$. The temporal landscape of $\varphi \imp \psi$ contains a roof over all time intervals within which $\varphi$ is contained in $\psi$, as shown in \cref{fig.temp_land_implication}.\footnote{In the following, we will often restrict ourselves to temporal landscapes on some arbitrary bounded interval, typically starting at 0, just for typographical convenience.} In the example shown in \cref{fig.temp_land_implication}, the landscape $\psi$ might appear having a confusing shape. The interpretation of such a landscape is that the predicate $\psi$ is only true over intervals of length at most three, the union of which is the solid line of height 3 above the axis.

Now that we have defined the logical connectives, we move on to the quantifiers $\exists(x:X)\ldotp P(x)$ and $\forall(x:X)\ldotp P(x)$. The simplest case is when these quantifiers range over a constant type $A$, which we may think of simply as a set.%
\footnote{In temporal type theory we can also quantify over non-constant behavior types, e.g.\ $\forall (x:X)\ldotp P(x)$, but this is a bit more technical. Since such quantification appears only in a single subsection (e.g.\ in \cref{eqn.agentpos}), we simply refer the reader to \cite{schultz2019temporal} for a definition.}
 Given a set $A$, a function $P\colon A\to\prop$ is a collection of $|A|$-many temporal landscapes. Taking their union defines the temporal landscape $\exists(a:A)\ldotp P(a)$, which will be a temporal landscape. Taking their intersection may not satisfy condition (b) of \cref{def.temp_landscape}, so we define $\forall(a:A)\ldotp P(a)$ to be the largest temporal landscape contained in this intersection.

Throughout this document, the reader will see the connectives $\wedge,\vee,\imp,\neg$, and the quantifiers $\exists$ and $\forall$. In each case, they refer to the operations on landscapes defined above.

\section{Predicates over grid worlds}\label{chap.examples}

In this section we give increasingly expressive examples of how to use temporal landscapes to describe the behavior of an agent moving in various types of environments. We start with a fairly standard model, used in the Artificial Intelligence (AI) literature, to describe motions of an agent over a discretized space. 

In a non-temporal situation, it is typical to represent an environment as a two- (or higher-)dimensional regular grid where each region of the space is a cell and cells overlap only on specified boundaries. Mathematically it is convenient to model this with an undirected graph $G = (V,E)$, where $V=\{0,\dots,N\}$ is a set of vertices, associated to subdivisions (or cells) of the environment, and $E\subseteq V\times V$ is the set of edges representing the fact that it is possible to move from one subdivision to another.

As we are working temporally, we replace sets with \emph{behavior types}, which may be thought of as time-varying sets. More precisely, a behavior type specifies, for every temporal landscape, a set of behaviors that could take place over those durations. These sets of behaviors are required to obey a certain compatibility condition, so that for example behaviors over long time intervals restrict to behaviors on shorter subintervals. 

\subsection{Modelling the environment: constant and non-constant behavior types}

\noindent\textbf{Static environments.} To simplify matters, let us first consider the case in which the environment does not vary in time. For this, we use constant behavior types: given a set $X$, the \emph{constant behavior type} on $X$, by abuse of notation written again simply as $X$, is the behavior type that for every temporal landscape simply specifies $X$ as its set of possible behaviors.

Suppose we want to say that our environment is modelled by the graph $(V,E)$, and that it does not change over time. To do this, we simply take $V$ and construct its constant behavior type $V$, and take $E$ as the constant subtype of $V \times V$ consisting precisely of the pairs $(v_1,v_2)$ in the set $E$. The constancy of the subtype $E$ says that the adjacency relation does not change: $v_1$ and $v_2$ either are adjacent or are not adjacent, independently of time.

To describe this fact logically in temporal type theory (TTT), we may write the formula
\[\forall(v_1,v_2:V)\ldotp(v_1,v_2)\in E\vee(v_1,v_2)\not\in E. 
\]
That said, in higher order logics like TTT, one typically exchanges subobjects for predicates, e.g.\ replacing $E\ss V\times V$ with $E\colon V\times V\to\prop$.\footnote{Recall that $\prop$ is the set of temporal landscapes; see \cref{def.temp_landscape}.} Then the statement would read 
\begin{equation}\label{eqn.decidable}
	\forall(v_1,v_2:V)\ldotp E(v_1,v_2)\vee\neg E(v_1,v_2).
\end{equation}
If we impose axiom \cref{eqn.decidable} then for any $v_1,v_2$, the landscape for $E(v_1,v_2)$ is either the \emph{always-true landscape} $\true$, or the \emph{always-false landscape} $\false$, depending on whether we want an edge $(v_1,v_2)$ or not. In this case we would say that $(V,E)$ forms a \emph{constant graph}.

\noindent\textbf{Dynamic environments.} It is also interesting, however, to decline to require that our environment obey axiom \cref{eqn.decidable}, and thus model environments in which the adjacency of cells changes over time. This makes sense in the autonomous setting if we imagine that sometimes a door is blocked or a secret passage is opened.

That said, for most situations, including all that follow, it is good enough to use a model of the environment where adjacency is \emph{symmetric}: if $v_1$ is connected to $v_2$, then $v_2$ is connected to $v_1$. Using the language of TTT, this means our environment obeys the axiom
\[\forall(v_1,v_2:V)\ldotp (v_1,v_2)\in E \iff (v_2,v_1) \in E\,.\]

It will also be convenient to work with the function $V \to (V\to\prop)$ given by \emph{currying} $E\colon V\times V\to\prop$. It sends a cell $v:V$ to the (time-varying) set $\{v':V\mid E(v,v')\}$ of cells adjacent to $v$. For our example, we want to consider the notion of \emph{neighbor}, by which we mean an adjacent cell, not including the cell itself. For each $v:V$, its neighbors are the subtype of $V$ defined by the formula
\[
	\pred{Nbr}(v)\coloneqq \{v': V\mid v'\neq v \wedge E(v,v')\}\,.
\]
It is worth noticing, once more, that we can interpret~$\pred{Nbr}$ in terms of temporal landscapes, by considering $\pred{Nbr}:V\to (V\to \prop)$. This means that $\pred{Nbr}(v)(v')$ is a truth value---i.e.\ a temporal landscape---that describes when a given~$v':V$ is a neighbor of~$v:V$. Of course, if we assume that $(V,E)$ is a constant graph, again $\pred{Nbr}(v)(v')$ will either be the $\true$ landscape or the $\false$ landscape, depending on whether $v$ and $v'$ are neighbors.

Working with TTT feels much like predicate logic and set theoretic constructions. In contrast with temporal logic, where one must get used to working with new logical operators such as `until' and `since'. For those who are trained in these languages, TTT provides a more easy to read language for reasoning about time. However, it is also true that TTT, being an intuitionistic logic, does introduce some ``complexities'', as the law of excluded middle (or equivalently, double negation elimination) needs not to hold. Temporal landscapes provide a graphical interpretation that helps reason about temporal logic statement and mitigate some of these complexities.

\subsection{Free/occupied cells: the negation operator}

We now expand our example by adding a predicate $\occ(v):V \to \prop$.
This predicate will be assumed to model the idea of a cell being \emph{occupied}: for each cell $v$, it specifies the set of time intervals over which $v$ is occupied. We will see that this predicate can capture situations that are more interesting than ``mere occupancy'', and that temporal landscapes provides a formal language to express such scenarios.

If a cell is not occupied, we will say that it is \emph{free}.
We further define $\free\coloneqq\neg\occ$; for each $v:V$, the landscape $\free(v)$ is the set of time intervals over which $v$ is free. As mentioned in \cref{sec.logical_manipulations}, double negation is not a trivial operation (the logic is constructive rather than Boolean). The predicate $\occ$ gives a good example of why this might be useful, i.e.\ why it makes sense that $\occ\overset{?}{=}\neg\neg\occ$ need not hold.
\label{rem.not_boolean}

Suppose we have agents $A$, $B$, and $C$, and predicates $\occ_A$, $\occ_B$, and $\occ_C$, which map a cell $v$ to the temporal landscape of intervals over which the respective agent is in $v$. Suppose that we wish to define $\occ$ to be the predicate describing the intervals over which at least one of the agents $A$, $B$, and $C$ is in $v$. Note that this is slightly ambiguous in English, but we will see that the negation operator in TTT allows us to easily distinguish between the two readings of this sentence as $\occ$ and $\neg\neg\occ$.

To do this, define $\occ$ to be the disjunction of these three predicates. For each $v$, $\occ(v)$ specifies the time intervals over which a single agent, whether it be $A$, $B$, or $C$, remains in the cell throughout. Then $\neg\neg\occ(v)$ specifies the time intervals over which there is always at least one agent in $v$, but agents are allowed to come and go.

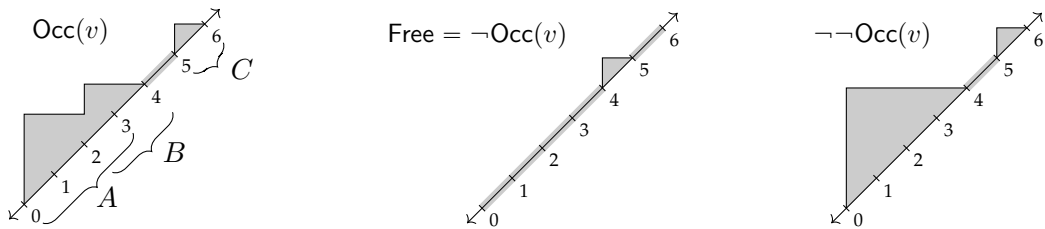
\begin{figure}
	\centering
	\subfloat[][]{
	\begin{tikzpicture}[scale=.4]
		\ArbitraryLandscape{(0,0) -- (0,3) -- (2,3) -- (2,4) -- (4,4)}		
		\False[gray!40]{4}{5}
		\True[gray!40]{5}{6}
		\LandscapeTimeAxis{0}{6}
		\node[below left, font=\small] at (current bounding box.north) {$\occ(v)$};
		\draw [decorate,decoration={brace,amplitude=5pt,mirror,raise=10pt}]
		(0,0) -- (3,3) node [black,midway,xshift=0.5cm, yshift=-0.9cm,label={$A$}] {};
		\draw [decorate,decoration={brace,amplitude=5pt,mirror,raise=15pt}]
		(2,2) -- (4,4) node [black,midway,xshift=0.8cm, yshift=-0.9cm,label={$B$}] {};
		\draw [decorate,decoration={brace,amplitude=5pt,mirror,raise=10pt}]
		(5,5) -- (6,6) node [black,midway,xshift=0.7cm, yshift=-0.8cm,label={$C$}] {};
	\end{tikzpicture}}
	\hspace*{.6in}
	\subfloat[][]{
	\begin{tikzpicture}[scale=.4]
		\False[gray!40]{0}{4}
		\True[gray!40]{4}{5}
		\False[gray!40]{5}{6}
		\LandscapeTimeAxis{0}{6}
		\node[below left, font=\small] at (current bounding box.north) {$\free(v)=\neg\occ(v)$};
	\end{tikzpicture}}
	\hspace*{.6in}
	\subfloat[][]{
	\begin{tikzpicture}[scale=.4]
		\True[gray!40]{0}{4}
		\False[gray!40]{4}{5}
		\True[gray!40]{5}{6}
		\LandscapeTimeAxis{0}{6}
		\node[below left, font=\small] at (current bounding box.north) {$\neg\neg\occ(v)$};
	\end{tikzpicture}}
	\caption{Temporal landscapes for: (a) $\occ(v)$, (b) $\free(v)\coloneqq\neg\occ(v)$, and (c) $\neg\free(v)=\neg\neg\occ(v)$.}\label{fig.occ_notocc_notnotocc}\vspace*{-0.6cm}
\end{figure}

More concretely, fix some cell $v$ and suppose that an agent~$A$ is in $v$ throughout the interval $[0,3]$, an agent $B$ is in $v$ throughout $[2,4]$, and another agent $C$ is in $v$ throughout $[5,6]$. Then the temporal landscapes for $\occ(v)$, for $\free(v)\coloneqq\neg\occ(v)$, and for $\neg\free(v)=\neg\neg\occ(v)$ are shown on the left, middle, and right, of \cref{fig.occ_notocc_notnotocc}.
While the middle and right-hand diagram fit the usual interpretation of ``when'' the room is free/occupied, they are derived from the left-hand diagram, which is more expressive. 

In particular, note that on the left-hand side diagram we have that $\occ(v)$ does not contain the time interval $[1.5,3.5]$, because this point falls in between the two roofs corresponding to $A$ and~$B$ occupying the cell. This might appear strange, since there is at least one agent in the cell throughout $[1.5,3.5]$, and thus we might expect $\occ(v)$ to contain this interval. However, $\occ$ expresses the more refined idea of those intervals over which there exists any specific agent occupying $v$: agent $A$ occupies $v$ on the interval $[1,3]$ and $B$ occupies $v$ on the interval $[2,4]$; $\occ(v)$ is the disjunction of $A$, $B$ and $C$.

As an example of where the extra expressivity of $\occ$ might be important, 
%
consider a simple situation where an emergency light is placed within the cell $v$ and where two consecutive ``blinks'' of the light would represent a dangerous situation. In the case the light is ON at $1.5$ and again at $3.5$, the temporal landscape for $\occ(v)$ would correctly capture the fact that such an alarm would be missed as there is not a single agent in the cell at those instances. Thus, unless $A$ and $B$ communicate about the status of the light, the notification of danger would be completely missed. Such communication---and memory /  recall in general---amounts to a strategy for persistently encoding intervallic facts into the present state.

\subsection{Objects in a room: quantifiers}

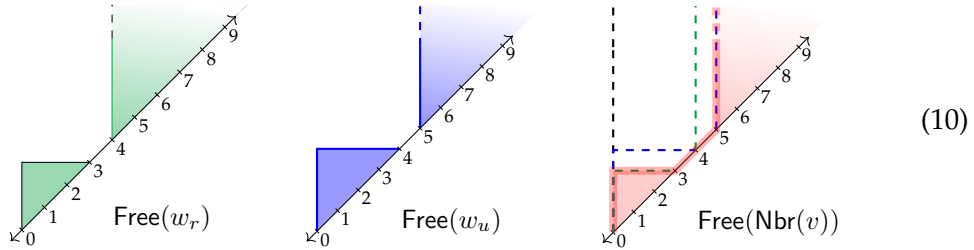
\begin{figure}[t!]
		\centering
		\subfloat[][]{\label{fig.grid_world_new}
		\begin{scaletikzpicturetowidth}{.21\textwidth}
		\begin{tikzpicture}[scale=\tikzscale, baseline=(current  bounding  box.center)]
			\fill[black!50] (-1.5,3.0) rectangle (0,4.5);
			\fill[black!50] (3,0) rectangle (4.5,1.5);
			\draw (-4.5,4.5) rectangle (4.5,-4.5);
			\foreach \x in {1,...,5}{
				\draw (-4.5+1.5*\x,-4.5) -- (-4.5+1.5*\x,4.5);
				\draw (-4.5, -4.5+1.5*\x) -- (4.5,-4.5+1.5*\x);
			}
			\draw node[fill, circle, minimum size=0.5cm, inner sep=0pt] at (-0.75,3.75) {};
			\draw[thick, ->] (-0.75,3.75) -- (-0.75,2.25);
			\draw node[fill, circle, minimum size=0.5cm, inner sep=0pt] at (3.75,0.75) {};
			\draw[thick, ->] (3.75,0.75) -- (2.25,0.75);
			\node at (-2.2,-0.8) {$v$};
			\node at (-3.7,-0.8) {$w_\ell$};
			\node[blue] at (-2.2,0.7) {$w_u$};
			\node[bluegreen] at (-0.8,-0.8) {$w_r$};
			\node at (-2.3,-2.3) {$w_d$};
		\end{tikzpicture}
	\end{scaletikzpicturetowidth}}
	\hfill
	\subfloat[][]{\label{fig.free_wr}
		\begin{scaletikzpicturetowidth}{.21\textwidth}
			\begin{tikzpicture}[scale=\tikzscale,baseline=(current  bounding  box.center)]
			\AlwaysTrueFrom[bluegreen][-8]{4}{10}
			\True[bluegreen!40]{0}{3}
			\LandscapeTimeAxis{0}{9}
			\node[right] at (3.6, 0.5) {\small $\free(w_r)$};
			\end{tikzpicture}
		\end{scaletikzpicturetowidth}}
	\hfill
	\subfloat[][]{\label{fig.free_wu}
		\begin{scaletikzpicturetowidth}{.21\textwidth}
			\begin{tikzpicture}[scale=\tikzscale,baseline=(current  bounding  box.center)]
			\fill[white, bottom color = blue!40, top color = blue!0, middle color=blue!20, shading angle=-8](5,5)--(5,11)--(11,11)--cycle;
			\fill[blue!40](0,0)--(0,4)--(4,4)--cycle;
			\draw[thick, blue] (0,0) -- (0,4) -- (4,4);
			\draw[thick, blue] (5,5) -- (5,9) edge[dashed] (5,11);
			\LandscapeTimeAxis{0}{9}
			\node[right] at (3.6, 0.5) {\small $\free(w_u)$};
			\end{tikzpicture}
		\end{scaletikzpicturetowidth}}
	\hfill
	\subfloat[][]{\label{fig.free_nbr_v}
		\begin{scaletikzpicturetowidth}{.21\textwidth}
			\begin{tikzpicture}[scale=\tikzscale,baseline=(current  bounding  box.center)]
			\fill[white, bottom color = red!20, top color = red!0, middle color=red!10, shading angle=-8](5,5)--(5,11)--(11,11)--cycle;
			\fill[red!20](0,0) -- (0,3 )-- (3,3) -- cycle;
			\LandscapeTimeAxis{0}{9}
			\draw[thick, dashed, black] (0,0) -- (0,9) edge[dashed] (0,11);
			\draw[thick, dashed, blue] (0,0) -- (0,4) -- (4,4);
			\draw[thick, dashed, blue] (5,5) -- (5,11);
			\draw[thick, dashed, bluegreen] (0,0) -- (0,3) -- (3,3);
			\draw[thick, dashed, bluegreen] (4,4) -- (4,11);
			\begin{scope}[transparency group, opacity=0.30]
			\draw[line width=3pt, red] (0,0) -- (0,3) -- (3,3) -- (5,5) -- (5,9) edge[dashed] (5,11);
			\end{scope}
			\node[right] at (3.6, 0.5) {\small $\free(\pred{Nbr}(v))$};
			\end{tikzpicture}
		\end{scaletikzpicturetowidth}}
	\caption{(a) Grid world. (b)-(d) Temporal landscapes related to the grid world model.}\label{fig.grid_world_dynamic}
\end{figure}

In this subsection we use TTT to describe when neighbors of a cell are free (unoccupied), despite possibly moving obstacles. We will consider two scenarios, both depicted in~\cref{fig.grid_world_new}.

This might be important, for example, for an autonomous vehicle, where one might want to know when it is immediately adjacent to an obstacle, and hence should be wary of a collision. In these scenarios the large black dots each represent an obstacle. In the first scenario the objects are stationary; in the second, they move in the direction shown by the arrow.

To begin, note that we can extend a predicate $\free$ over any subtype $N\colon V\to\prop$ as follows:\footnote{Note that $\free(N)$ is (constructively) equivalent to $\forall(v:V)\ldotp \occ(v)\imp\neg N(v)$.}
\begin{align}\label{eqn.extension}
	\free(N) \coloneqq \forall (v:V)\ldotp N(v)\imp\free(v)\,.
\end{align}
This predicate describes the intervals over which all $v \in N$ are free. In particular, we will be interested in $\free(\pred{Nbr}(v))$ for a cell $v$, which tells us when every neighbor of $v$ is free, or equivalently, when none of $v$'s neighbors is occupied.


\noindent\textbf{Static objects.} 
Assume that the objects, represented by the two large black dots, are static (i.e.\ forget the arrows in \cref{fig.grid_world_new} for now). Consider the cell $v$ indicated in \cref{fig.grid_world_new}. In the case that the black objects are static, one sees that the predicate $\free(\pred{Nbr}(v))$ is the always-true landscape $\true$, since the configuration of ``free'' cells does not change over time. In particular, if we were to draw the temporal landscapes $\free(w)$ for each $w:\pred{Nbr}(v)$, each one would be the always-true landscape, and so their conjunction.

\noindent\textbf{Dynamic objects.} Next we consider a situation in which the black dots represent moving objects. In this scenario, the two objects move in the indicated directions (one downwards and the other leftwards) at a rate of one cell per unit time. When they reach a cell adjacent to the boundary of the domain, they remain there forever after.

In this case there is an equality of predicates $\free(w_\ell)=\free(w_d)$; both correspond to the always-true landscape, because these two cells are never occupied by either of the moving obstacles. For the predicate $\free(w_r)$, however we note that the cell will be occupied for one time unit, between 3 and 4. Thus the cell $w_r$ is free for the interval $[0,3]$, and of course for any subinterval of it, such as $[1, 1.46]$. The cell $w_r$ is also free for any interval $[4,b]$, as long as $4<b$; the temporal landscape $\free(w_r)$ is shown in \cref{fig.free_wr}. Similar reasoning applied to $\free(w_u)$ yields the temporal landscape shown in \cref{fig.free_wu}.

The temporal landscape for $\free(\pred{Nbr}(v))$ is the conjunction of these temporal landscapes, i.e.\ all of $v$'s neighbors must be free, as shown in \cref{fig.free_nbr_v}. As described in \cref{sec.logical_manipulations} the resulting temporal landscape is going to be the ``minimum'' landscape---red boldface in \cref{fig.free_nbr_v}---of the four landscapes shown with dashed lines $\free(w_\ell)$ and $\free(w_d)$ in black, $\free(w_r)$ in green, and $\free(w_u)$ in blue.

Looking at the red landscape on the right of \cref{fig.grid_world_dynamic} we immediately see that the neighborhood of~$v$, namely~$\pred{Nbr}(v)$, is not free in the interval $[3,5]$.

\begin{figure}[t!]
	\centering
	\subfloat[][]{\label{fig.grid_world_room}	
	\begin{scaletikzpicturetowidth}{0.21\textwidth}
		\begin{tikzpicture}[scale=\tikzscale,baseline=(current  bounding  box.center)]
			\fill[gray!45] (-1.5,-4.5) rectangle (0,-1.5);
			\fill[gray!45] (-1.5,0) rectangle (0,4.5);
			\fill[gray!45] (3,-4.5) rectangle (4.5,4.5);
			\fill[gray!45] (0,3) rectangle (3,4.5);
			\fill[gray!45] (0,-3) rectangle (3,-4.5);
			\fill[blue!12] (0,-3) rectangle (3,3);
			\draw (-4.5,4.5) rectangle (4.5,-4.5);
			\foreach \x in {1,...,5}{
				\draw (-4.5+1.5*\x,-4.5) -- (-4.5+1.5*\x,4.5);
				\draw (-4.5, -4.5+1.5*\x) -- (4.5,-4.5+1.5*\x);
			}
			\node at (0,-0.5) {$\rightarrow$};
			\node at (0,-1) {$\leftarrow$};
			\node at (0.75,0) {$\uparrow$};
			\node at (0.75,1.5) {$\uparrow$};
			\node at (1.5,2.25) {$\rightarrow$};
			\node at (2.25,0) {$\downarrow$};
			\node at (2.25,1.5) {$\downarrow$};
			\node at (1.5,-0.75) {$\leftarrow$};
			\draw node[pin=150:$a$, draw=blue, fill=blue, thick, circle, minimum size=0.3cm, inner sep=0pt] at (-0.75,-0.75) {};
		\end{tikzpicture}
	\end{scaletikzpicturetowidth}}
	\hfill
	\subfloat[][]{\label{fig.time_in_room}
		\begin{scaletikzpicturetowidth}{.21\textwidth}
			\begin{tikzpicture}[scale=\tikzscale,baseline=(current  bounding  box.center)]	
				\False[blue!40][0.5]{0}{3}
				\True[blue!40][0.5]{3}{6}
				\False[blue!40][0.5]{6}{7}
				\True[blue!40][0.5]{7}{15}
				\UnionAlwaysTrue[red!40][0.5]{0}{15}{3}
				\LandscapeTimeAxis{0}{15}
			\end{tikzpicture}
	\end{scaletikzpicturetowidth}}
	\hfill
	\subfloat[][]{\label{fig.implication}
		\begin{scaletikzpicturetowidth}{.21\textwidth}
			\begin{tikzpicture}[scale=\tikzscale,baseline=(current  bounding  box.center)]
				\ArbitraryLandscape{(0,0) -- (0,10) -- (7,10) -- (12,15) -- (15,15)}
				\draw[dashed] (7,7) --(7,10);
				\draw[dashed] (7,10) -- (10,10);
				\LandscapeTimeAxis{0}{15}
			\end{tikzpicture}
	\end{scaletikzpicturetowidth}}
	\caption{(a) Grid world with a room whose walls are shown in gray and interior in light blue. (b) Temporal landscape for $\forall(a:A)\ldotp (\pred{Pos}(a)\ss R)$ (blue) and $\exists(s:\rr)\ldotp \timebetw{s}{s+3}$) (red). (c) Temporal landscape of the implication \eqref{eq:time_in_room}.}\label{fig.agent_in_room_for_tau_sec}
\end{figure}


\subsection{Max dwell time in a room: implication} 

Let $A$ denote the type of agents' IDs and let $\pred{Pos}\colon A\to(V\to\prop)$ denote the predicate that an agent $a:A$ is at a vertex $v:V$. We also define a room~$R$ to be a subset of vertices, $R\colon V\to\prop$.\footnote{Note that we have not said $A$ and $R$ are constant over time: agents might come into service or be decommissioned, and rooms might be constructed, demolished, or expanded over time.} Then for an agent~$a:A$, one may write $\pred{Pos}(a)\ss R$ to denote the proposition $\forall(v:V)\ldotp\pred{Pos}(a)(v)\imp R(v)$, namely the agent~$a$ is in the room $R$. The situation is shown in \cref{fig.grid_world_room}, where we indicate an agent by a blue dot and shade in lighter blue the cells forming a room $R$.

The wall---cells that are always occupied---are depicted in gray. Arrows depict possible trajectories that agent~$a$ can take to move within the room~$R$ and then exit.

Suppose we want to express the proposition that an agent stays in a room~$R$ for at most~$\tau$ units of time before it must exit the room. To model this, for some $\tau:\rr_{\geq0}$, we can use the predicate
\begin{align}
\forall(a:A)\ldotp (\pred{Pos}(a)\ss R)\imp\exists(s:\rr)\ldotp \timebetw{s}{s+\tau}\,,\label{eq:time_in_room}
\end{align}
which says that given an agent~$a$, as long as~$a$'s position remains in room~$R$, there is some start time~$s$ such that the clock remains between~$s$ and~$s+\tau$. 


Let us consider an example. Let $\tau=3$ and suppose the agent is not in the room during the intervals $[0,3]$ and $[6,7]$, but is in the room during $[3,6]$ and $[7,15]$. The landscapes for the left and right hand side of~\eqref{eq:time_in_room}, namely $\pred{Pos}\ss R$ and $\exists(s:\rr)\ldotp s<t<s+\tau$ respectively, are shown in blue and red in~\cref{fig.time_in_room}.

Note the temporal landscape of the right hand side of~\eqref{eq:time_in_room} (red), it is an ``always true'' capped at 3 time units. This is because given a time~$t$ there always exists a real value~$s$ with $t\in[s,s+\tau]$, and such predicate is true for intervals $[t_1,t_2]$ of length $t_2-t_1\leq 3$.

The temporal landscape for the entire predicate~\eqref{eq:time_in_room}, is shown in \cref{fig.implication}. Note that in the interval [0,10] it is always true that the agent is within the room for at most 3 time units, however in [7,15] it is never true that the agent is within the room for at most 3 time units. Indeed, all we can say is that on intervals of length at most 3, the agent is clearly is in the room for at most 3 time units, however this is not true for longer time intervals.

\subsection{Regions and occupancy}

Let us now consider a modified grid world, where instead of constructing a uniform spatial partitioning of the environment we leverage what we, as humans, would argue is a reasonable ``semantic'' subdivision of the space. Take for example a building, such a partitioning would be based on more abstract concepts than cells, such as ``rooms'', ``corridors'', ``foyers'' etc. To ground the discussion, consider the scenario depicted in \cref{fig.building}. This picture depicts an agent (shown with a red dot) traversing a continuous environment, following a trajectory shown by a dashed line. The environment has been semantically subdivided into different regions (rooms). 

In order to model the agent moving through the building, we begin by saying what a trajectory is. We normalize the building to be the square $S\coloneqq[0,6]^2\ss\rr^2$. Then define the set of all possible time-parametrized trajectories through the building to be:
\[
  \mathcal{X} \coloneqq  \left\{(x_1,x_2):\wt{\rr}\times\wt{\rr}\mid 0\leq x_1\leq 6\text{ and } 0\leq x_2\leq 6\right\}\,.
\]
For example, in \cref{fig.building}, we depict a time-parametrized trajectory over an interval $(t_0,t_{11}$), where $t_0=0$, $t_1=1$, $t_2=2$, etc.\ making the distance traveled per unit time non-uniform along the trajectory.

It is worth noticing that in earlier examples, we considered a discretized space, whereas now we are considering a continuous space $S$, and behaviors defined over such space. 

\begin{figure}[t!]
	\centering
	\subfloat[][]{\label{fig.building}
		\begin{scaletikzpicturetowidth}{.32\textwidth}
			\begin{tikzpicture}[scale=\tikzscale,baseline=(current  bounding  box.center)]
				\draw[white] (-6,-6) rectangle (7,7){};
				\fill[black!50] plot coordinates{(3,1.5) (3,3) (0,3) (0,4.5) (4.5,4.5) (4.5,1.5)};
				\fill[black!50] (3,0) rectangle (4.5,1.5);
				\fill[orange!20] (-4.5,1.5) rectangle (-1.5,4.5); 
				\fill[blue!20] (-1.5,0) rectangle (3,3) ;
				\fill[yellow!20] (-4.5,-3)rectangle (1.5,0); 
				\fill[yellow!20] (-4.5,0)rectangle (-1.5,1.5);
				\fill[green!20] (-4.5,-3) rectangle (4.5,-4.5);
				\fill[green!20] (1.5,0) rectangle (4.5,-3);
				\fill[red!20] (-1.5,3) rectangle (0,4.5);
				\draw[line width=2pt, black!50] (-1.5,1.5) -- (-4.5,1.5) -- (-4.5,4.5) -- (-1.5,4.5) -- (-1.5,3);
				\draw[line width=2pt, black!50] (1.5,0) -- (-1.5,0) -- (-1.5,1.5) -- (-4.5,1.5) -- (-4.5,-3) -- (1.5,-3) -- (1.5,-1.5);
				\draw[line width=2pt, black!50] plot coordinates{(-4.5,-4.5) (4.5,-4.5) (4.5,0)};
				\draw[dashed] plot[smooth, tension=.7, mark=x, mark options={color=red,solid}, mark size=7pt] coordinates {(-1,4.5) (-0.6,2.8) (-2.2,2.3) (2.3,1.1) (1.5,-0.7) (-0.8,-0.7) (-3.5,0.4) (-3.3,-2.5) (2.8,-0.8) (2.7,-3.6) (0,-4) (-4.4,-3.8)};
				\node[label={[label distance=-7pt]45:{\tiny $t_0$}}] at (-1,4.5) {};
				\node[label={[label distance=-7pt]-45:{\tiny $t_1$}}] at (-0.6,2.8) {};
				\node[label={[label distance=-5pt]-90:{\tiny $t_2$}}] at (-2.2,2.3) {};
				\node[label={[label distance=-7pt]45:{\tiny $t_3$}}] at (2.3,1.1) {};
				\node[label={[label distance=-7pt]155:{\tiny $t_4$}}] at (1.5,-0.7) {};
				\node[label={[label distance=-7pt]45:{\tiny $t_5$}}] at (-0.8,-0.7) {};
				\node[label={[label distance=-7pt]-45:{\tiny $t_6$}}] at (-3.5,0.4) {};
				\node[label={[label distance=-7pt]65:{\tiny $t_7$}}] at (-3.3,-2.5) {};
				\node[label={[label distance=-7pt]45:{\tiny $t_8$}}] at (2.8,-0.8) {};
				\node[label={[label distance=-7pt]-45:{\tiny $t_{9}$}}] at (2.7,-3.6) {};
				\node[label={[label distance=-7pt]45:{\tiny $t_{10}$}}] at (0,-4) {};
				\node[label={[label distance=-7pt]45:{\tiny $t_{11}$}}] at (-4.4,-3.8) {};
				\node at (-2,-1) {\scriptsize $\mathsf{Room}_B$};
				\node at (-3,4) {\scriptsize $\mathsf{Room}_A$};
				\draw[dotted] (-1.5,3) -- (0,3);
				\draw[dotted] (-1.5,4.5) -- (0,4.5);
				\draw[dotted] (-1.5,3) -- (-1.5,1.5);
				\draw[dotted] (1.5,0) -- (3,0);
				\draw[dotted] (1.5,0) -- (1.5,-1.5);
				\draw[dotted] (-4.5,-4.5) -- (-4.5,-3);
				\node[pin={[pin distance=0.45cm, pin edge={thick, black}]40:{\scriptsize $\mathsf{Entrance}$}}] at (-0.6,4) {};
				\node at (0.5,0.8) {\scriptsize $\mathsf{Lobby}$};
				\node at (-1.2,-3.4) {\scriptsize $\mathsf{Corridor}$};
				\draw node [draw=blue, fill=blue, circle, minimum size=0.3cm, inner sep=0pt] at (-0.8,2) {};
				\node[label={[label distance=-1pt]90:0}] at (-4.5,4.5) {};
				\node[label={[label distance=-1pt]90:6}] at (4.5,4.5) {};
				\node[label={[label distance=-1pt]0:6}] at (4.5,4.5) {};
				\node[label={[label distance=-1pt]0:0}] at (4.5,-4.5) {};
			\end{tikzpicture}
	\end{scaletikzpicturetowidth}}
	\hspace*{2cm}
	\subfloat[][]{\label{fig.agent_in_room}
		\begin{scaletikzpicturetowidth}{.32\textwidth}
			\begin{tikzpicture}[scale=\tikzscale,baseline=(current  bounding  box.center)]
				\draw[white] (-2,-2) rectangle (13,13){};
				\True[red!20][0.5]{0}{11/5*0.9}
				\True[blue!20][0.5]{11/5*0.7}{11/5*1.4}
				\True[orange!20][0.5]{11/5*1.2}{11/5*2.4}
				\True[blue!20][0.5]{11/5*2.2}{11/5*4.6}
				\True[green!20][0.5]{11/5*4.4}{11/5*5.1}
				\LandscapeTimeAxisNoTicks{0}{11}
				\draw (0+.1,0-.1) to[pos=-3, font=\tiny] node {$t_0$} (0-.1, 0+.1);
				\draw (11/5*1+.1,11/5*1-.1) to[pos=-3, font=\tiny] node {$t_1$} (11/5*1-.1,11/5*1+.1);
				\draw (11/5*2+.1,11/5*2-.1) to[pos=-3, font=\tiny] node {$t_2$} (11/5*2-.1,11/5*2+.1);
				\draw (11/5*4+.1,11/5*4-.1) to[pos=-3, font=\tiny] node {$t_3$} (11/5*4-.1,11/5*4+.1);
				\draw (11/5*5+.1,11/5*5-.1) to[pos=-3, font=\tiny] node {$t_4$} (11/5*5-.1,11/5*5+.1);
				\draw[dotted] (11/5*0.8-11/5*0.8,11/5*0.8+11/5*0.8) -- (11/5*0.8,11/5*0.8);
				\draw[dotted] (11/5*1.3-11/5*0.8,11/5*1.3+11/5*0.8) -- (11/5*1.3,11/5*1.3);
				\draw[dotted] (11/5*2.3-11/5*0.8,11/5*2.3+11/5*0.8) -- (11/5*2.3,11/5*2.3);
				\draw[dotted] (11/5*4.5-11/5*0.8,11/5*4.5+11/5*0.8) -- (11/5*4.5,11/5*4.5);
				\draw[dotted] (11/5*5-11/5*0.8,11/5*5+11/5*0.8) -- (11/5*5,11/5*5);
			\end{tikzpicture}
	\end{scaletikzpicturetowidth}}
	\caption{(a) Layout of a building with the trajectory (dashed line) followed by an agent (blue dot). Specific locations where the agent is at time instances $t_0, t_1, \dots$ are shown with red crosses. (b) Temporal landscape corresponding to the agent being in a room. The colors correspond to the rooms in (a).}
\end{figure}



As in the discrete case, suppose we are given a predicate $\occ\colon S\to\prop$ that models the subset of $S$ (possibly changing in time) in which the agent cannot be. For example in Figure~\ref{fig.building} we show some gray regions, which are intended to be walls, and hence always occupied/non-traversable. Thus if $s:S$ is in a wall, we put $\occ(s)\coloneqq\true$. Again as in the discrete case, we use this to define a predicate $\free\colon\mathcal{X} \to \prop$, by $\free(x)\coloneqq\forall(s:S)\ldotp\occ(s)\imp\neg(x=s)$.

The agent is moving through the building, but to be more realistic we could imagine that the agent occupies space larger than a point, and that different parts of the agent move at slightly different speeds. Hence the agent consists of several different trajectories, all of which are close to one another, say within a distance of $\gamma:\rr$. We define $\pred{close}\colon\mathcal{X}\times\mathcal{X}\to\prop$ using the Euclidean norm $\pred{close}(x_1,x_2)\coloneqq \|x_1-x_2\|_2\leq\gamma$. We also put an upper bound on the speed of the agent, say $\velbnd:\rr$, and define the agent's possible positions as the following behavior type:
\begin{equation}\label{eqn.agentpos}
	\type{AgentPos} \coloneqq 
	\left\{\pred{p}:\mathcal{X}\to \prop \;\middle|\; 
	\parbox{3.0in}{\raggedright$
		\forall (x_1, x_2:\mathcal{X})\ldotp\big((\pred{p}(x_1)\wedge \pred{p}(x_2)) \imp \pred{close}(x_1,x_2)\big) \;\wedge$\\$ \forall(x:\mathcal{X})\ldotp \pred{p}(x) \imp \big(\free(x) \wedge -\velbnd\leq\dot{x} \leq\velbnd\big)$}
	\right\}\,.
\end{equation}
Here the bound $-\velbnd\leq\dot{x}\leq\velbnd$ is the temporal landscape consisting of those intervals $[t_1,t_2]$ over which $\|x(t_1')-x(t_2')\|_2<(t_2'-t_1')\velbnd$ holds for all $t_1<t_1'<t_2'<t_2$. For more on derivatives in temporal type theory, see \cite[Section 7.3]{schultz2019temporal}.

Let us consider the constant type $R\coloneqq\{\mathsf{Room_A},\allowbreak\mathsf{Room_B},\allowbreak\mathsf{Entrance},\allowbreak\mathsf{Lobby},\allowbreak\mathsf{Corridor}\}$ representing the rooms as shown in Figure~\ref{fig.building}. Suppose we have a predicate $\pred{Room}\colon R\to(\mathcal{X}\to\prop)$, indicating the landscape on which a trajectory stays within a room. As before, for any moving agent $a:A$ with trajectories $\pos(a):\type{AgentPos}$, let us denote with $\pos(a) \subseteq r$ the predicate   
$\forall(x:\mathcal{X})\ldotp \pos(a)(x) \imp \pred{Room}(r)(x)$, which says that agent $a$ is in a room~$r$ if all of the trajectories that make up $a$ are in~$r$.

The temporal landscape for the proposition $\pos(a)\ss r$ for when a single agent $a$ is in a room $r$ is a roof. Thus the temporal landscape of $\pred{AgentInARoom} := \exists(r:R)\ldotp \pos(a)\subseteq r$ is obtained by taking the union---namely the max---of these roofs as shown in \cref{fig.agent_in_room}.

Note that because of the agent footprint there are intervals where the agent can be in two rooms and since the agent is always in some room and thus $\neg\neg\pred{AgentInARoom}$ is the always-true landscape.

\subsection{Landmarks and maps: ``slanted'' temporal landscapes}

So far, for most of the examples---except for that in \cref{fig.agent_in_room_for_tau_sec}---all the temporal landscapes have consisted of a finite union of roofs. One thus wonders when a ``slanted'' temporal landscape would be relevant in an application and what it would represent.

\begin{figure}[t!]
	\centering
	\subfloat[][]{\label{fig.slanted_tl_env}
	\begin{scaletikzpicturetowidth}{0.7\textwidth}
		\begin{tikzpicture}[scale=\tikzscale,baseline=(current  bounding  box.center)]
		\draw[thick]  plot coordinates {(1, 0) (2.5, 0) (2.5,0.6) (4.5,0.6)(4.5,0)(5.2,0)};
		\draw[thick]  plot coordinates {(1,-0.6) (2.5,-0.6) (2.5,-1.2) (3.9,-1.2)(4.5,-1.2)(4.5,-0.6)(5.2,-0.6)};
		\draw[thick, dashed, red, opacity=0.5, ->]  plot coordinates {(1,-0.3) (5.1,-0.3)};
		\draw[draw=white, fill=blue, opacity=0.2]  (1.6,-0.3) ellipse (0.6 and 0.6);
		\draw[draw=white, fill=blue, opacity=0.2]  (3.2,-0.3) ellipse (0.6 and 0.6);
		\foreach \a in {5, 10,..., 360} {
			\draw[opacity=0.2] {(1.6,-0.3)-- +(\a:0.6)};
			\draw[opacity=0.2] {(3.2,-0.3)-- +(\a:0.6)};}
		\draw node [draw=red, fill=red, circle, minimum size=0.4cm, inner sep=0pt,pin={[pin distance=1.5cm]135:{1}}] at (1.6,-0.3) {};
		\draw node [draw=red, fill=red, circle, minimum size=0.4cm, inner sep=0pt,pin={[pin distance=1.8cm]45:{2}}] at (3.2,-0.3) {};
		\draw[fill=black]  (2.9,0) rectangle (3,0.3);
		\draw[fill=black]  (3.4,0) rectangle (3.3,0.3);
		\draw[fill=black]  (3.3,-0.6) rectangle (3.4,-0.9);
		\draw[fill=black]  (2.9,-0.6) rectangle (3,-0.9);
		\draw[fill=black]  (4,-0.6) rectangle (4.1,-0.9);
		\draw[fill=black]  (4,0) rectangle (4.1,0.3);	
		\foreach \i in {0,0.1,...,0.2} {
			\draw node [draw=blue, fill=blue, circle, minimum size=0.10cm, inner sep=0pt] at (2.9+\i,0) {};
		}
		\foreach \i in {0,0.05,...,0.1} {
			\draw node [draw=blue, fill=blue, circle, minimum size=0.10cm, inner sep=0pt] at (2.9+\i,0) {};
			\draw node [draw=blue, fill=blue, circle, minimum size=0.10cm, inner sep=0pt] at (3.3+\i,0) {};
			\draw node [draw=blue, fill=blue, circle, minimum size=0.10cm, inner sep=0pt] at (4+\i,0) {};
			\draw node [draw=blue, fill=blue, circle, minimum size=0.10cm, inner sep=0pt] at (2.9+\i,-0.6) {};
			\draw node [draw=blue, fill=blue, circle, minimum size=0.10cm, inner sep=0pt] at (3.3+\i,-0.6) {};
			\draw node [draw=blue, fill=blue, circle, minimum size=0.10cm, inner sep=0pt] at (4+\i,-0.6) {};
		}
		
		\foreach \i in {0,0.05,...,0.7} {
			\draw node [draw=blue, fill=blue, circle, minimum size=0.1cm, inner sep=0pt] at (4.5+\i,0) {};
			\draw node [draw=blue, fill=blue, circle, minimum size=0.10cm, inner sep=0pt] at (4.5+\i,-0.6) {};
		}
		\foreach \i in {0,0.05,...,1.4} {
			\draw node [draw=blue, fill=blue, circle, minimum size=0.10cm, inner sep=0pt] at (1.1+\i,0) {};
			\draw node [draw=blue, fill=blue, circle, minimum size=0.10cm, inner sep=0pt] at (1.1+\i,-0.6) {};
		}
		\foreach \i in {0,0.05,0.1,0.2,0.3} {
			\draw node [draw=blue, fill=blue, circle, minimum size=0.10cm, inner sep=0pt] at (2.5,\i) {};
			\draw node [draw=blue, fill=blue, circle, minimum size=0.10cm, inner sep=0pt] at (2.5,-\i-0.6) {};
			\draw node [draw=blue, fill=blue, circle, minimum size=0.10cm, inner sep=0pt] at (4.5,\i) {};
			\draw node [draw=blue, fill=blue, circle, minimum size=0.10cm, inner sep=0pt] at (4.5,-\i-0.6) {};
			\draw node [draw=blue, fill=blue, circle, minimum size=0.10cm, inner sep=0pt] at (2.9,\i) {};
			\draw node [draw=blue, fill=blue, circle, minimum size=0.10cm, inner sep=0pt] at (3,\i) {};
			\draw node [draw=blue, fill=blue, circle, minimum size=0.10cm, inner sep=0pt] at (3.3,\i) {};
			\draw node [draw=blue, fill=blue, circle, minimum size=0.10cm, inner sep=0pt] at (3.4,\i) {};
			\draw node [draw=blue, fill=blue, circle, minimum size=0.10cm, inner sep=0pt] at (4,\i) {};
			\draw node [draw=blue, fill=blue, circle, minimum size=0.10cm, inner sep=0pt] at (4.1,\i) {};
			\draw node [draw=blue, fill=blue, circle, minimum size=0.10cm, inner sep=0pt] at (2.9,-\i-0.6) {};
			\draw node [draw=blue, fill=blue, circle, minimum size=0.10cm, inner sep=0pt] at (3,-\i-0.6) {};
			\draw node [draw=blue, fill=blue, circle, minimum size=0.10cm, inner sep=0pt] at (3.3,-\i-0.6) {};
			\draw node [draw=blue, fill=blue, circle, minimum size=0.10cm, inner sep=0pt] at (3.4,-\i-0.6) {};
			\draw node [draw=blue, fill=blue, circle, minimum size=0.10cm, inner sep=0pt] at (4,-\i-0.6) {};
			\draw node [draw=blue, fill=blue, circle, minimum size=0.10cm, inner sep=0pt] at (4.1,-\i-0.6) {};  
		}
		\end{tikzpicture}
	\end{scaletikzpicturetowidth}}
	\hspace*{2cm}
	\subfloat[][]{\label{fig.slanted_temp_land}
		\begin{scaletikzpicturetowidth}{0.22\textwidth}
		\begin{tikzpicture}[scale=\tikzscale,baseline=(current  bounding  box.center)]
			\draw[draw=white, bottom color = gray!40, top color = gray!0, middle color=gray!20, shading angle=315]  plot coordinates {(0,0) (0,3.6)(1.8,5.4)(2.8,7.4) (3.4,9.4)(6.2,10)(8.4,12)(12,12)};
			\draw  plot coordinates {(0,0) (0,3.6)(1.8,5.4)(2.8,7.4) (3.4,9.4)(6.2,10)(8.4,12)};
			\LandscapeTimeAxis{0}{10}
		\end{tikzpicture}
	\end{scaletikzpicturetowidth}}
	\caption{(a) LIDAR returns (blue dots) as an agent (red dot) travels along a corridor. (1) and (2) show two specific locations of the agent as it travels left to right. (b) Example of a slanted temporal landscape.}
\end{figure}

Consider the scenario as shown in~\cref{fig.slanted_tl_env}. An agent (red dot) travels, at constant velocity, within an indoor environment along the red dashed path. The agent is equipped with a range limited sensor, such as a LIDAR (blue disk) which emits a set of discrete laser beams. For each laser beam the LIDAR gets a return whenever a laser beam hits a surface. The measurement is the (possibly noisy) location of the surface along each beam (blue dots). We have denoted with (1) and (2) two specific locations along the path. We depict with small blue dots a possible set of sensor measurements---samples---along walls and columns (black rectangles).

Further suppose that the agent is equipped with a fixed amount of onboard memory, so that not all the samples can be stored. When the buffer used to store samples is full, past samples will be deleted to make space for new ones. Identify each sample with an integer $i:\nn$ and define the predicate $\pred{SampleInMem}(i)$ that will be true over $[t_1,t_2]$ as long as the sample~$i$ is in the memory of the agent.

The temporal landscape for $\pred{SampleInMem}(i)$ is clearly a roof over the interval $[t_1,t_2]$ where $t_1$ is the instance when the sample~$i$ was first stored in memory and $t_2$ is the instance when it was overwritten by a new sample (for $t_2 = +\infty$ then there is enough memory so that no overwriting occurs).

The following predicate will have a ``slanted'' temporal landscape 
\[
\pred{SamplesInMem} = \bigvee_i \pred{SampleInMem}(i)\,.
\]
For example, it might look like the temporal landscape in~\cref{fig.slanted_temp_land}. Initially, in the corridor, the number of samples is high and the memory will be fully allocated. As new samples are obtained, old ones will be overwritten. Assuming a constant velocity and number of samples per unit of time, we have a constant overwriting so that a sample is in memory only over a constant size interval: thus the landscape will be parallel to the time line. As the agent enters a part of the environment that has fewer surfaces, the number of samples  per unit time decreases, and thus samples will persist in memory over longer and longer periods of time, especially given that the environment becomes sparser as the agent moves left to right. Once the agent starts sensing the beginning of the right-most corridor, the number of samples starts to quickly increase, the persistence of a sample in memory decreases and the maximum persistence in the onboard memory is reached (the landscape is parallel to the time line again).

It is worth mentioning that an idea related to temporal landscapes has been applied to other problems in robotics and in particular to monitor and diagnose perception systems~\cite{PA-DIS-LC:21}. In this context, the concept of temporal diagnostic graphs was introduced. 

\section{Conclusion}\label{chap.conclusion}

In this paper, we have attempted to give an intuitive introduction to the logic of temporal type theory in terms of temporal landscapes. On the one hand, these are just collections of time intervals over which a proposition may be true. On the other, they can be drawn as Lipschitz functions and hence visualized. They form a logical system, where all of the connectives and quantifiers are defined by operations on these Lipschitz functions.

After introducing these landscapes, we discussed a series of examples from the domain of autonomous agents. These became fairly complex, e.g.\ considering an agent's position not just as a point but as a collection of points, each moving with bounded speed, avoiding possibly moving obstacles, and storing recent LIDAR measurements in a small-capacity memory that is constantly being overwritten. These examples point to the great expressivity of temporal type theory.

In practice, TTT can serve as a sort of big tent, where calculations from model checkers or ODE solvers can be embedded. While infinite in nature, we explained how temporal landscapes can be finitely approximated. We thus hope to have shown how TTT can be used to specify and guide algorithmic developments in autonomous systems and any modeling environment in which time is an issue. 

\printbibliography

\end{document}